\documentclass[oneside,a4paper,12pt,reqno]{amsart}
\usepackage{upref,amsfonts,amsxtra,amssymb}

\newcommand{\Z}{\mathbb{Z}}
\newcommand{\T}{\mathbb{T}}
\newcommand{\R}{\mathbb{R}}
\newcommand{\g}{\Gamma}
\newcommand{\tg}{\mathbb{T}^{\Gamma}}
\newtheorem{theorem}{Theorem}[section]
\newtheorem{lemma}[theorem]{Lemma}

\newtheorem{prop}[theorem]{Proposition}

\title{Ergodicity of algebraic actions of nilpotent groups}
\author{Siddhartha Bhattacharya}
\address{School of Mathhematics, Tata Institute of Fundamental
  Research, Mumbai 400005, India}
\email{siddhart@math.tifr.res.in}
\date{}
\begin{document}
\maketitle
\begin{abstract}
An algebraic $\Gamma$-action is an action of a countable group $\Gamma$
on a compact abelian group $X$ by continuous automorphisms of $X$. We
prove that any expansive algebraic action of a finitely generated
nilpotent group $\Gamma$ on a connected group $X$ is ergodic. We also
show that this result does not hold for actions of polycyclic groups.   
\end{abstract}
\section{Introduction}
Let $\Gamma$ be a countable group. An {\it algebraic $\Gamma$-action\/} 
is an action $\alpha$ 
of $\Gamma$ on a compact metrizable abelian group $X$
by continuous automorphisms of $X$. Any such action
induces a $\Gamma$-action on ${\widehat X}$, the dual of $X$, by
automorphisms of ${\widehat X}$. Hence ${\widehat X}$ can be viewed as
a $\Z[\Gamma]$-module, where $\Z[\Gamma]$ is the integral group ring
of $\Gamma$. Using duality theory one can show that this
identification gives rise to a bijective correspondence between
algebraic $\g$-actions and modules over $\Z[\g]$. 
The action $\alpha$ is called {\it  Noetherian\/} if
${\widehat X}$ is Noetherian as a $\Z[\Gamma]$-module. Equivalently, 
$\alpha$ is Noetherian if any decreasing sequence
$$ X = X_0 \supset X_1\supset X_2\supset\cdots $$
of closed $\Gamma$-invariant subgroups stabilizes.

It is easy to see that any algebraic $\g$-action $(X,\alpha)$
preserves the Haar measure on $X$. Hence $\alpha$ can
also be viewed as a measure preserving dynamical system. 
For $d\ge 1$, ergodic algebraic $\Z^d$-actions were characterized in
\cite{KS}. For non-abelian $\g$ this problem has been studied under
the additional assumption that the action $\alpha$ is a principal
action, i.e., the dual module of $\alpha$ is of the form 
$\Z[\g]/\Z[\g]f$ for some $f\in \Z[\g]$. It is known that for a large
class of groups, that includes free groups with more than one generators
as well as all finitely generated amenable groups that are not virtually
cyclic, such actions are always ergodic (\cite{H}, \cite{LPS}). However, when $\g$ is
non-abelian, the ergodicity properties of non-principal algebraic
$\g$-actions are not well understood, even when $\g$ is the discrete
Heisenberg group (see \cite[Problem 3.5]{LS}).
 
In this paper we prove the following result :
\begin{theorem}\label{main}
Let $(X,\alpha)$ be an  expansive algebraic action
of a finitely generated nilpotent group $\g$.
 Then the
$\sigma$-algebra of $\alpha$-invariant sets is finite. If $X$ is
also connected, then $\alpha$ is ergodic. 
\end{theorem}
In \cite{CL} Chung and Li obtained several characterizations of
expansiveness for algebraic actions of countable groups. Using their
result we construct a large class of actions where Theorem
\ref{main} can be applied.
Chung and Li also raised the question whether algebraic
quotients of expansive actions of polycyclic-by-finite groups are
always expansive (\cite[Conjecture 3.6]{CL}). We 
show that this is not true
in general. More specifically, we construct a
polycyclic group $\g\subset GL(3,\Z)$ such that the natural action of
$\g$ on $\T^3$ is expansive but it admits a non-trivial quotient on which
the induced $\Gamma$-action is trivial. This proves that Theorem
\ref{main} does not hold for actions of polycyclic groups.

 \section{Ergodicity and Expansiveness}

Our proof of Theorem \ref{main} relies on the observation that under
the stated assumptions the action $\alpha$ exhibits a certain form of 
cohomological rigidity.
 A similar technique, but involving a
different form of cohomological rigidity,
 was used in \cite{B2} to study periodic
points of algebraic actions. Let us first recall a few related definitions. 
Let $\alpha$ be an action of a countable group $\Gamma$ on a compact
metrizable abelian group $X$ by continuous automorphisms of $X$. A 
{\it 1-cocycle \/} of $\alpha$ is a map $c : \Gamma\rightarrow X$
 that satisfies the equation 
$$ c(\gamma \gamma^{'}) = c(\gamma) + \alpha(\gamma )(c(\gamma^{'}))$$
for all $\gamma, \gamma^{'}\in \Gamma$.
We will denote the collection of all 1-cocycles by $C(\alpha)$.
If $X^{\Gamma}$ denote the compact abelian group of all functions from 
$\Gamma$ to $X$, equipped with the product topology and pointwise
addition, then it is easy to see that $C(\alpha)\subset X^{\Gamma}$ is 
a closed subgroup.
For any $x\in X$, the map $c_x : \Gamma\rightarrow X$ defined by 
$c_x(\gamma) = \alpha(\gamma)(x) - x$ is a 1-cocycle. 
A cocycle $c$ is said to be a {\it coboundary \/} if
$c = c_x$ for some $x\in X$. Since the map $x\mapsto c_x$ is a
continuous homomorphism from $X$ to $C(\alpha)$, it follows that 
$B(\alpha)$, the collection of all 
coboundaries, is a closed subgroup of $C(\alpha)$. 
Two 1-cocycles 
$c_1$ and $c_2$ are said to be {\it cohomologous } if $c_1 - c_2\in B(\alpha)$. 
The quotient group $C(\alpha) /B(\alpha)$ is denoted
by $H^{1}(\alpha)$ and is called the {\it first cohomolgy group}
 of the action $\alpha$.
\begin{lemma}\label{ce}
Let $(X,\alpha)$ be an algebraic action of a countable group $\Gamma$,
and let $K\subset X$ be a closed $\Gamma$-invariant subgroup. If
$\beta$ denotes the induced $\Gamma$-action on $X/K$ then $H^{1}(\alpha)$
is finite whenever both $H^{1}(\beta)$ and $H^{1}(\alpha|_{K})$ are finite. 
\end{lemma}

{\it Proof. } Let $\pi$ denote the quotient homomorphism from $X$ to
$X/K$. We define a homomorphism $\theta$ from $C(\alpha)$
to $C(\beta)$ by $\theta(c) = \pi\circ c$. 
As $\pi\circ c_x = c_{\pi(x)}$, $\theta$ induces a homomorphism
$\theta^{*}$ from $H^{1}(\alpha)$ to $H^{1}(\beta)$. 
 Suppose $c$ is an element of
$C(\alpha)$ such that the corresponding element of $H^{1}(\alpha)$
lies in ${\rm Ker}(\theta^{*})$. Then $\pi\circ c = c_y$ for some $y\in
X/K$. We pick $q\in X$ with $\pi(q) = y$. As $\pi\circ(c - c_q)
= 0$, $c$ is comohomologous to a 1-cocycle taking values in
$K$. Since $H^{1}(\alpha|_{K})$ is finite, we deduce that that ${\rm
  Ker}(\theta^{*})$ is finite. Since ${\rm Image }(\theta^{*})\subset
H^{1}(\beta)$ is also finite, this proves the given assertion. $\Box$

\medskip
For any algebraic action $(X, \alpha)$ of a countable group $\Gamma$
let $F(\alpha)\subset X$ denote the set of points that are fixed by
every element of $\alpha(\Gamma)$. Clearly $F(\alpha)$
is a closed $\g$-invariant subgroup of $X$. 
\begin{lemma}\label{cq}
Let $(X,\alpha)$ be an algebraic action of a countable group $\Gamma$,
and let $K\subset X$ be a closed $\Gamma$-invariant subgroup such that
the induced $\Gamma$-action on $X/K$ has infinitely many fixed points.
Then either $F(\alpha)$ is infinite or $H^{1}(\alpha|_{K})$ is
infinite. 
\end{lemma}
{\it Proof.} Suppose $H^{1}(\alpha|_{K})$ is finite. 
Let $\beta$ denote the quotient action of $\g$ on $X/K$, and let
$\pi : X\rightarrow X/K$ denote the projection map.
If $Y\subset X$ denotes the
subgroup $\pi^{-1}(F(\beta))$, then $Y$ is $\Gamma$-invariant, 
$\Gamma$ acts trivially on $Y/K$, and $Y/K$ is infinite. 
We define a homomorphism 
$\theta$ from $Y$ to $C(\alpha|_{K})$ by
$\theta(y)(\gamma) = \alpha(\gamma)(y) - y$. 
 Since 
$H^{1}(\alpha|_{K})$ is finite, $B(\alpha|_{K})$ is a finite index
subgroup of $C(\alpha|_{K})$. Let $Z\subset Y$ denote the
$\Gamma$-invariant subgroup $\theta^{-1}(B(\alpha|_{K}))$. Then
$Z$ contains $K$, and since 
$Y/Z$ is finite, $Z/K$ is infinite.
We pick an arbitrary $z\in Z$ and choose $p\in K$ such that 
$\theta(z) = c_{p}$. Then $\alpha(\gamma)(z - p ) = z - p$ 
for all $\gamma\in\Gamma$, i.e., $z- p$ is fixed by $\alpha$. This
shows that every coset of $K$ in $Z$ contains an element of
$F(\alpha)$. As $Z/K$ is infinite, this shows that $F(\alpha)$ is
infinite. 
$\Box$ 

\medskip
Now we consider the case when $\g$ is a finitely 
generated nilpotent group. In this case there exists a sequence of subgroups
$$ \{ e\} = \Gamma_0 \subset \Gamma_1\subset \cdots \subset \Gamma_n =
\Gamma$$
such that for all $i$, $\Gamma_i$ is normal in $\Gamma_{i+1}$ and 
$\Gamma_{i+1}/\Gamma_i$ lies in the center of $\Gamma_{n}/\Gamma_i$. 
Since $\Gamma_{i+1}/\Gamma_i$ is a finitely generated
abelian group for each $i$, refining the sequence if necessary, we may
assume that $\Gamma_{i+1}/\Gamma_i$ is cyclic. Let $l(\Gamma)$ denote
the smallest possible length for such a series. We note that for any
finitely generated nilpotent group $\g$  there exists 
$\gamma_0$ in the center of $\g$ such that $l(\g / \Z\gamma_0) <
l(\g)$, where $\Z\gamma_0$ is the cyclic subgroup generated by $\gamma_0$.
 
\begin{theorem}\label{nilp} Let $\g$ be a finitely generated
  nilpotent group, and let $(X,\alpha)$ be a Noetherian action of $\g$
  such that $F(\alpha)$ is finite. Then $H^{1}(\alpha)$ is finite.
\end{theorem}
{\it Proof.} We will use induction on $l(\g)$. 
If $l(\g) = 0$
then $\g = \{ e\}$. 
Since $F(\alpha)$ is finite, this shows that $X$ is finite
 and the above assertion is trivially 
true. Suppose $l(\g)\ge 1$ and $(X,\alpha)$ is a Noetherian action
of $\g$ such that $F(\alpha)$ is finite but $H^{1}(\alpha)$ is
infinite. Let ${\mathcal A}$ denote the collection of all closed
$\Gamma$-invariant subgroups $L$ such that $H^{1}(\alpha|_{L})$ is
infinite. 
Since $\alpha$ is Noetherian, ${\mathcal A}$ has a minimal element $K$. 
We choose $\gamma_0$  
in the center of $\g$ such that
$l(\g / \Z\gamma_0) < l(\g)$. Let $L\subset X$ denote the
$\g$-invariant subgroup $(\alpha(\gamma_0) - I)(K)$. We claim that 
$K = L$. Suppose this is not the case. Then $H^{1}(\alpha|_{L})$ is finite.
Let $\beta$ denote the quotient action of $\g$ on $K/L$. Clearly
$\beta$ is a Noetherian action.
As $F(\alpha|_{K})$ is
finite, applying Lemma \ref{cq} we obtain that $F(\beta)$ is finite. 
Since $\beta(\gamma_0) = I$, for any 1-cocycle $r\in C(\beta)$ and
$\gamma\in\g$ we have,
$$ r(\gamma_0\gamma) = r(\gamma_0) + \beta(\gamma_0)r(\gamma)
= r(\gamma_0) + r(\gamma).$$
As $\gamma_0\gamma = \gamma\gamma_0$, we also have,
$  r(\gamma_0\gamma) = 
 r(\gamma) + \beta(\gamma)r(\gamma_0)$.
This shows that $r(\gamma_0)\subset F(\beta)$. Let $\theta :
C(\beta)\rightarrow F(\beta)$ denote the group homomorphism defined by
$\theta(r) = r(\gamma_0)$. As ${\rm Image}(\theta)\subset F(\beta)$ is
finite, we deduce that ${\rm Ker}(\theta)$ is a finite-index subgroup
of $C(\beta)$. Since $\beta(\gamma_0) = I$, $\beta$ induces an
automorphism action $\beta^{*}$ of $\g/\Z\gamma_0$ on $K/L$. It is
easy to see that any 
$s\in {\rm Ker}(\theta)$ induces a 1-cocycle 
${\overline{s}}\in C(\beta^{*})$. Moreover, $s$ is a coboundary if and
only if ${\overline{s}}$ is a coboundary. As $l(\g/\Z\gamma_0) <
l(\g)$, by the induction hypothesis $H^{1}(\beta^{*})$ is
finite. Hence  ${\rm Ker}(\theta)$ contains a finite index subgroup
that consists of coboundaries. As ${\rm Ker}(\theta)$ is a finite
index subgroup of $C(\beta)$, this shows that  $H^{1}(\beta)$ is
finite. Applying Lemma \ref{cq} we conclude that $H^{1}(\alpha|_{K})$
is finite. This contradicts the defining property of $K$ and proves
the claim.

Now let $H\subset K$ denote the subgroup defined by
$$ H = \{ x\in K : \alpha(\gamma_0)(x) = x\}.$$
As $(\alpha(\gamma_0) - I)(K) = K \ne \{ 0\}$, $H$ is a proper closed
subgroup of $K$. Since $\gamma_0$ lies in the center of $\Gamma$, $H$
is also $\Gamma$-invariant. From the minimality of $K$ we deduce that 
$H^{1}(\alpha|_{H})$ is finite. Now let $c$ be an arbitrary
$1$-cocycle of $\alpha|_{K}$. We choose $p\in K$ such that 
$(\alpha(\gamma_0) - I)(p) = c(\gamma_0)$ and define $c^{'}\in
C(\alpha|_{K})$ by $c^{'} = c - c_{p}$. Since $c^{'}(\gamma_0) = 0$,
it follows that for any $\gamma\in\Gamma$,
$$c^{'}(\gamma\gamma_0) = c^{'}(\gamma) +
\alpha(\gamma)(c^{'}(\gamma_0)) = c^{'}(\gamma).$$
On the other hand, we also have
$$c^{'}(\gamma_0\gamma) = c^{'}(\gamma_0) +
\alpha(\gamma_0)(c^{'}(\gamma)) = \alpha(\gamma_0)(c^{'}(\gamma)) .$$
Since $\gamma_0$ commutes with every element, this shows that 
the image of $c^{'}(\gamma)$ is contained in $H$. Hence every
1-cocycle of $\alpha|_{K}$ is cohomologous to a 1-cocycle taking
values in $H$. As $H^{1}(\alpha|_{H})$ is finite we deduce that 
 $H^{1}(\alpha|_{K})$ is also finite. This contradicts the
 defining property of $K$ and proves the induction step. $\Box$

\medskip
{\it Proof of Theorem \ref{main} :} Since $(X,\alpha)$ is expansive,
${\widehat{X}}$ is finitely generated as a $\Z[\g]$-module. As $\g$ is
a finitely generated nilpotent group, $\Z[\g]$ is a Noetherian ring. 
Hence $\widehat{X}$ is a Noetherian $\Z[\g]$-module, i.e., 
the action $\alpha$ is Noetherian.
Let $M\subset {\widehat{X}}$
denote the submodule consisting of all characters $\chi$ such that
the $\g$-orbit of $\chi$ under the dual action $\widehat{\alpha}$  is
finite. Let $\pi$ denote the projection map from $X$ to
$\widehat{M}$. Then $Y = {\rm Ker}(\pi)$ is a closed
$\Gamma$-invariant subgroup and the restriction of $\alpha$ to $Y$ is
ergodic. Since $\alpha$ is Noetherian, $M$ is finitely generated. Let
$\{ \chi_1, \ldots ,\chi_n\}\subset M$ be a generating set, and for $i =
1, \ldots ,n$; let $\Lambda_i\subset\Gamma$ be the finite-index
subgroup that fixes $\chi_i$. Then every element of $M$ is fixed by
the finite-index subgroup $\Lambda = \cap \Lambda_i$, i.e., $\Lambda$
acts trivially on ${\widehat{M}} = X/Y$. Let $\beta$
denote the restriction of $\alpha$ to $\Lambda$. From the
expansiveness of $\alpha$ it follows that $F(\beta)\subset
X$ and $F(\beta|_{Y})\subset Y$ are finite sets. 
By Theorem \ref{nilp}, $H^{1}(\beta|_{Y})$ is also finite. 
As every element of $X/Y$ is fixed by $\beta$, applying Lemma \ref{cq}
we see that $X/Y$ is finite. Now the ergodicity of $\alpha|_{Y}$
implies that the invariant $\sigma$-algebra of $\alpha$ is finite. If
$X$ is connected then it does not admit non-trivial finite-index
subgroups. Hence $Y = X$ , i.e., $\alpha$ is ergodic. $\Box$     
\section{Examples}
Let us first describe a class of actions where Theorem \ref{main} can be
applied. For any countable group $\g$ let $\tg$ denote the compact
abelian group of all functions from $\g$ to $\T$, equipped with
pointwise addition and the product topology. The shift action $\alpha$
of $\g$ on $\tg$ is defined by
$\alpha(\gamma)(x)(\gamma_1) = x(\gamma^{-1}\gamma_1)$. For $f = 
\sum c_{\gamma}\delta_{\gamma}\in \Z[\g]$ and $x\in \tg$ we define
$$ f(x) = \sum c_{\gamma }x(\gamma).$$
For a left ideal $I\subset  \Z[\g]$  we define 
$X(I) = \{ x\in\tg : f(x) = 0\ \forall f\in I\}$. Conversely, if 
$X\subset \tg$ is a closed shift-invariant subgroup we set 
$I_{X} = \{ f\in \Z[\g] : f(x) = 0 \ \forall x\in X\}$. From the
duality theory it follows that both $I\mapsto X(I)$ and $X\mapsto I_X$
are bijective correspondences between left ideals of $\Z[\g]$ and 
the closed shift-invariant subgroups of $\tg$. 

For a countable group $\g$ let $l^{1}(\g)$ denote the Banach algebra
consisting of all
functions $h$ from $\g$ to ${\mathbb{C}}$ with $\sum_{\gamma}|h_{\gamma}|
< \infty$. An element $f = \sum c_{\gamma}\delta_{\gamma}\in \Z[\g]$
is called {\it lopsided} if
there exists $\gamma_0$ such that $|c_{\gamma_0}| > \sum_{\gamma\ne
  \gamma_0}|c_{\gamma}|$. One can show that any lopsided $f\in\Z[\g]$
is invertible in $l^{1}(\g)$, and any left ideal $I\subset \Z[\g]$ 
containing an  invertible element of $l^{1}(\g)$ also contains a 
lopsided element of $\Z[\g]$
(\cite[Proposition 5.4]{LS}).  From \cite[Theorem 3.1]{CL} it follows
that for a left ideal $I\subset \Z[\g]$, the shift action of $\g$ on
$X(I)$ is expansive if and only if $I$ contains an invertible element
of $l^{1}(\g)$. In particular, for any ideal $I\subset \Z[\g]$
containing a lopsided element of $\Z[\g]$ the shift action of $\g$ 
on $X(I)\subset \tg$ is expansive, and Theorem \ref{main} can be
applied. 

For a countable group $\g$, and a
left ideal $I\subset \Z[\g]$ we define $I^{*}\subset \Z[\g]$ by
$$ I^{*} = \{ f\in \Z[\g] : nf\in I \ \mbox{for some}\ n\ge 1\}.$$ 
If $Y$ denotes the connected component of $X(I)$ that contains $0$
then from duality theory it follows that $I^{*} = I_{Y}$. Hence
if $J\subset \Z[\g]$ is a left ideal containing a lopsided element 
then by Theorem \ref{main} 
the algebraic $\g$-action corresponding the module $\Z[\g]/J^{*}$ is
ergodic.
  
\medskip
We will now show that Theorem \ref{main} does not hold for actions of
polycyclic groups. We will use the following result which 
is a special case of
\cite[Theorem A]{B1}. 
\begin{prop}
Let $n\ge 1$, $\Gamma\subset GL(n,\Z)$ be a subgroup, 
and $\alpha$ be the
natural action of $\Gamma$ on $\T^n\cong \R^n/\Z^n$. Then $\alpha$ is
expansive if and only if for every non-zero $p\in\R^n$ the $\Gamma$-orbit of
$p$ is unbounded. 
\end{prop}

\medskip
{\it Example :} 
We define $A\in GL(2,\Z)$ and a subgroup 
$\Gamma\subset GL(3,\mathbb{Z})$
 by
$$ A =   
\left[
\begin{array}{cc}
 2 & 1  \\
 1 & 1  \\
\end{array}
\right],\  
\Gamma =  \left\{  
\left(
\begin{array}{cc}
 A^n & b  \\
 0 & 1  \\
\end{array}
\right) 
\ \mid\  n\in\Z, b\in\Z^2\ 
 \right\}.  $$
Let $\Gamma_0\subset \Gamma$ denote the subgroup defined by
$$\Gamma_0 =  \left\{  
\left(
\begin{array}{cc}
 I & b  \\
 0 & 1  \\
\end{array}
\right) 
\ \mid\  b\in\Z^2\ 
 \right\}.  $$
It is easy to see that $\Gamma_0\subset\Gamma$ is normal and $\g/\Gamma_0\cong
\Z$. Since $\Gamma_0$ is isomorphic with $\Z^2$ this shows that $\g$ is
polycyclic.
Let $\alpha$ denote the natural action of $\g$ on $\T^3$. We claim
that $\alpha$ is expansive but not ergodic. 
By the above proposition, to prove expansiveness of $\alpha$ it is enough to
show that every non-trivial $\g$-orbit in $\R^3$ is unbounded. Let 
$p = (x,y,z)\in \R^3$ be a point with bounded $\g$-orbit. 
Then the $\Gamma_0$-orbit of $p$ is also a bounded. On
the other hand,
$$\Gamma_{0}p = \{ (x+mz, y+nz,z) : m,n\in\Z\}.$$
This shows that $z= 0$. Now we define $\Gamma_1\subset \g$ by
$$ \Gamma_1 =  \left\{  
\left(
\begin{array}{cc}
 A^n & 0  \\
 0 & 1  \\
\end{array}
\right) 
\ \mid\  n\in\Z\ 
 \right\}.  $$
Since $z = 0$, we see that $\Gamma_{1}p = \{ (A^{n}(x,y),0) :
n\in\Z\}$. Since $\Gamma_{1}p$ is a bounded subset of $\R^3$ and both
the eigenvalues of $A$ have modulus different from $1$, we deduce that 
$x = y = 0$. This proves that $\alpha$ is expansive.
We define $K = \{ (x,y,z)\in\T^3 : z = 0\} \subset \T^3$. It is easy
to see that $K$ is a $\g$-invariant closed subgroup and  $\g$ acts
trivially on $X/K$.  
Hence $\alpha$ is not ergodic.

\medskip
{\it Remark .} If $\alpha$ and $K\subset \T^3$ are as above, then from
the expansiveness of $\alpha$ and triviality of the quotient action on
$X/K$ it follows that the action $\alpha|_{K}$ does not have
pseudo-orbit tracing property. This shows that expansive algebraic
actions of polycyclic groups may have non-expansive algebraic
quotients and need not satisfy pseudo-orbit tracing property, giving
negative answers to questions raised by Chung and Li (\cite[Conjecture
3.6]{CL}) and
Meyerovitch (\cite[Question 3.11]{M}).

\end{document}